\documentclass[11pt,a4paper]{article}
\usepackage{amsmath}
\usepackage{graphicx}
\usepackage{amsfonts}
\usepackage{amssymb}
\newtheorem{theorem}{Theorem}

\newtheorem{corollary}[theorem]{Corollary}

\newtheorem{lemma}[theorem]{Lemma}

\begin{document}

\title{Poincar\'{e}'s theorem for the modular group of real Riemann surfaces}
\author{Antonio F. Costa
\and Sergey M. Natanzon\\{\small Both authors partially supported by MTM2005-01637}\\{\small Second author partially supported by }\\{\small NWO-RFBR 047.011.2004.026 (RFBR 05-02-89000-NWO\_a),}\\{\small NSh-1972.2003.1, RFBR 04-01-00762}\\}
\date{}
\maketitle

\textbf{Abstract}. Let $Mod_{g}$ denote the modular group of (closed and
orientable) surfaces $S$ of genus $g$. Each element $[h]\in Mod_{g}$ induces a
symplectic automorphism $H([h])$ of $H_{1}(S,\mathbb{Z})$. Poincar\'{e} showed
that $H:Mod_{g}\rightarrow Sp(2g,\mathbb{Z})$ is an epimorphism. A real
Riemann surface is a Riemann surface $S$ together with an anticonformal
involution $\sigma$. Let $(S,\sigma)$ be a real Riemann surface,
$Homeo_{g}^{\sigma}$ be the group of orientation preserving homeomorphisms of
$S$ such that $h\circ\sigma=\sigma\circ h$ and $Homeo_{g,0}^{\sigma}$ be the
subgroup of $Homeo_{g}^{\sigma}$ consisting of those isotopic to the identity.
The group $Mod_{g}^{\sigma}=$ $Homeo_{g}^{\sigma}/Homeo_{g,0}^{\sigma}$ plays
the r\^{o}le of the modular group in the theory of real Riemann surfaces. In
this work we describe the image by $H$ of $Mod_{g}^{\sigma}$. Such image
depends on the topological type of the involution $\sigma$.

\section{Introduction}

The modular group, $Mod_{g},$ of genus $g$ closed (compact and without
boundary) orientable surfaces (or mapping class group) is the group of
orientation preserving autohomeomorphisms of a surface $S$ of genus $g$ modulo
the homeomorphisms isotopic to the identity. Each element $[h]\in Mod_{g}$
induces in $H_{1}(S,\mathbb{Z})$ a symplectic automorphism $H([h])$, in such a
way that $H:Mod_{g}\rightarrow Sp(2g,\mathbb{Z})$ is a homomorphism. In the
celebrated paper \cite{P} by Poincar\'{e}, one of the principal results is
that $H$ is an epimorphism. This important relation between the modular group
of complex curves and the symplectic group has been used, for instance, in the
classification of abelian actions on surfaces, see \cite{CN} and \cite{E}.

The moduli space of complex algebraic curves is a central object in
mathematics and recently in mathematical physics. The moduli space is the
quotient of the Teichm\"{u}ller space by the action of the modular group. By
this fact, the modular group $Mod_{g}$ plays an important r\^{o}le in the
study of the geometry and topology of the moduli space of Riemann surfaces of
complex algebraic curves, see \cite{N} and \cite{MS}.

A real Riemann surface $(S,\sigma)$ is a Riemann surface $S$ with an
anticonformal involution $\sigma$. Real Riemann surfaces correspond to real
algebraic curves and the moduli space of real Riemann surfaces is the moduli
of real algebraic curves or ``real moduli space''. The study of the moduli
space of real algebraic curves is somehow similar to the study of the complex
case but it presents some difficulties of its own. The ``real moduli space''
has importance in its own right and has applications in many areas (see
\cite{D}). The modular group of a real Riemann surface (real modular group) is
the group defined as follows. Let $(S,\sigma)$ be a real Riemann surface,
$Homeo_{g}^{\sigma}$ be the group of orientation preserving homeomorphisms of
$S$ such that $h\circ\sigma=\sigma\circ h$ and $Homeo_{g,0}^{\sigma}$ be the
subgroup of $Homeo_{g}^{\sigma}$ consisting of homeomorphisms isotopic to the
identity. The real modular group $Mod_{g}^{\sigma}$ is the quotient
$Homeo_{g}^{\sigma}/Homeo_{g,0}^{\sigma}$. The real moduli for real Riemann
surfaces with the same topological type than $(S,\sigma)$ is the quotient of a
contractible space by the action of $Mod_{g}^{\sigma}$ (see \cite{N}). Our aim
is to describe the automorphisms of $H_{1}(S,\mathbb{Z})$ defined by
$Mod_{g}^{\sigma}$, i. e. to establish the analog of Poincar\'{e} theorem for
the real modular group.

Of course the involution $\sigma$ induces a nontrivial involution $\sigma_{H}$
on the homology group $H_{1}(S,\mathbb{Z})$. The image of the real modular
group by the homomorphism $H,$ defined by the action on homology, lies in
$Sp(2g,\mathbb{Z})\cap\{h:h\circ\sigma_{H}=\sigma_{H}\circ h\}=Aut^{\sigma
_{H}}(H_{1}(S,\mathbb{Z})).$ A naive conjecture would be that the image by $H$
of $Mod_{g}^{\sigma}$ is exactly $Aut^{\sigma_{H}}(H_{1}(S,\mathbb{Z}))$. This
turn out to be the case when the fixed point set of $\sigma$ is empty (Theorem
7). When the fixed point set of $\sigma$ is nonempty and consists of a
collection of disjoint simple closed curves, however, the subspace of
$H_{1}(S,\mathbb{Z})$ spanned by the fixed curves must be preserved. This last
condition is not sufficient and somewhat more technical conditions must be
imposed. In general $H(Mod_{g}^{\sigma})$ depends on the topological type of
$\sigma$. We now give precise statements.

Orientation reversing involutions $\sigma$ of $S$ up topological equivalence
are classified by an integer $\pm k$: the \textit{topological type}. The set
$Fix(\sigma)$ is called the \textit{real part} of $\sigma$ and it consist of
$k$ disjoint simple closed curves, called \textit{ovals}. The set
$Fix(\sigma)$ can be \textit{separating} if $S-Fix(\sigma)$ is disconnected
and we say that $\sigma$ has topological type $+k$. In this case
$S/\left\langle \sigma\right\rangle $ is orientable, $0<k\leq g+1$ and
$k\equiv g+1\operatorname{mod}2$. If the set $Fix(\sigma)$ is
\textit{non-separating}, that is to say $S-Fix(\sigma)$ is connected, then
$\sigma$ has topological type $-k$. Now $S/\left\langle \sigma\right\rangle $
is nonorientable and $0\leq k\leq g$ (Harnack theorem, see\textbf{\ }%
\cite{BEGG}).

As before, let $H:Mod_{g}^{\sigma}\rightarrow Sp(2g,\mathbb{Z})$ denote the
homomorphism given by the action on $H_{1}(S,\mathbb{Z})$ of the elements of
$Mod_{g}^{\sigma}$. Let $(.,.)$ be the intersection form in $H_{1}%
(S,\mathbb{Z})$. The principal results are as follows:

\bigskip

\textit{Case 1} (Theorem 7). If the topological type of $(S,\sigma)$ is $0$ then

\begin{center}
$H(Mod_{g}^{\sigma})=\{h\in Sp(2g,\mathbb{Z}):h\circ\sigma_{H}=\sigma_{H}\circ h\}.$
\end{center}

\bigskip

\textit{Case 2 }(Theorem 8). If the topological type of $(S,\sigma)$ is $-k$
then there is a symplectic basis $\left\{  X_{i},Y_{i}\right\}  _{i=1,...,g}$
in $H_{1}(S,\mathbb{Z}),$ i. e. $(X_{i},X_{j})=(Y_{i},Y_{j})=0$, $(X_{i}%
,Y_{j})=\delta_{ij}$, such that the first $k$ elements $X_{i}$ are represented
by the $k$ ovals (of $Fix(\sigma)$) and $h\in H(Mod_{g}^{\sigma})$ iff

1. $h\in Sp(2g,\mathbb{Z}),$

2. $h\circ\sigma_{H}=\sigma_{H}\circ h,$

3. $h\{X_{1},...,X_{k}\}=\{\varepsilon_{1}X_{1},...,\varepsilon_{k}X_{k}\},$
where $\varepsilon_{i}\in\{-1,1\}$, $i=1,...,k$.

4. if $i\in\{k+1,...,g\},$ $(h(X_{i}),X_{j})=0,$ for every $0\leq j\leq g$ and
$(h(X_{i}),Y_{j})$ is even for every $1\leq j\leq k$.

\bigskip

\textit{Case 3 }(Theorem 10). If the topological type of $(S,\sigma)$ is $+k$,
let the connected components of $S-Fix(\sigma)$ be denoted by $S_{1}$ and
$S_{2}$. Let $H_{i}\leq H_{1}(S,\mathbb{Z})$ be spanned by cycles represented
by closed curves in $S_{i}$, $i=1,2$. Let $X_{1},...,X_{k}$ be the $k$
elements of $H_{1}(S,\mathbb{Z})$ represented by the ovals (of $Fix(\sigma)$).

Then $h\in H(Mod_{g}^{\sigma})$ iff

1. $h\in Sp(2g,\mathbb{Z})$,

2. $h\circ\sigma_{H}=\sigma_{H}\circ h,$

3. $h\{X_{1},...,X_{k}\}=\{\varepsilon X_{1},...,\varepsilon X_{k}\}$ where
$\varepsilon=\pm1$ and $h(H_{i})\subset H_{i}$ if $\varepsilon=1$ and
$h(H_{i})\subset H_{j},$ $i\neq j$, if $\varepsilon=-1.$

\bigskip

\section{Preliminaries on crystallographic\ groups}

We use 2-dimensional crystallographic groups as a tool for our arguments. The
results are attainable using only topology arguments (generators of the
mapping class group of surfaces) but the difficulties are comparable.

Let $\mathbb{U}^{2}=\mathbb{S}^{2}$, $\mathbb{E}^{2}$ and $\mathbb{H}^{2}$ be
the sphere, euclidean plane and the hyperbolic plane respectively and let
$\mathcal{G}$ be the group of isometries of $\mathbb{U}^{2}.$ Let
$\mathcal{G}^{+}$ be the index two subgroup of $\mathcal{G}$ consisting of the
orientation preserving transformations of $\mathbb{U}^{2}$.

A crystallographic group, is a discrete subgroup $\Gamma$ of the group
$\mathcal{G}$ with compact quotient space. If the group $\Gamma$ contains
orientation reversing transformations then $\Gamma^{+}=\Gamma\cap
\mathcal{G}^{+}$ is an index two subgroup of $\Gamma$ and $\Gamma^{+}$ is the
canonical orientation preserving crystallographic subgroup of $\Gamma$.

If $\Gamma$ is such a group then its algebraic structure and the geometrical
structure of the quotient orbifold $\mathbb{U}^{2}/\Gamma$ is given by the
signature (see \cite{BEGG}):

\begin{center}
$s({\Gamma})=(h;\pm;[m_{1},...,m_{r}];\{(n_{11},...,n_{1s_{1}}),...,(n_{k1}%
,...,n_{ks_{k}})\}).$ \ \ \ \ \ \ \ \ \ \ \ \ \ \ \ \ (*)
\end{center}

The orbit space $\mathbb{U}^{2}/\Gamma$ is an orbifold with underlying surface
of genus $h$, having $r$ cone points and $k$ boundary components, each with
$s_{j}\geq0,$ $j=1,...,k,$ corner points. The signs $^{\prime\prime}%
+^{\prime\prime}$ and $^{\prime\prime}-^{\prime\prime}$ correspond to
orientable and non-orientable underlying surfaces respectively. The integers
$m_{i}$ are called the proper periods of ${\Gamma}$ and they are the orders of
the cone points of $\mathbb{U}^{2}/\Gamma$. The brackets $(n_{i1}%
,...,n_{is_{i}})$ are the period cycles of ${\Gamma}$ and the integers
$n_{ij}$ are the link periods of ${\Gamma}$; they are the orders of the corner
points of $\mathbb{U}^{2}/\Gamma$. The group $\Gamma$ is isomorphic to the
orbifold fundamental group $\pi_{1}O(\mathbb{U}^{2}/\Gamma)$.

A group $\Gamma$ with signature (*) has a \textit{canonical presentation} with generators:

$x_{1},...,x_{r},e_{1},...,,e_{k}$, $c_{ij},1\leq i\leq k,0\leq j\leq s_{i}$ and

$a_{1},b_{1},...,a_{h},b_{h}$ if $\mathbb{U}^{2}/\Gamma$ is orientable or

$d_{1},...,d_{h}$ otherwise,

\hspace{-0.5cm}and relators:\newline $x_{i}^{m_{i}},i=1,...,r$, $c_{ij}%
^{2},(c_{ij-1}c_{ij})^{n_{ij}},c_{i0}e_{i}^{-1}c_{is_{i}}e_{i}$,
$i=1,...,k,j=0,...,s_{i}$ and

$x_{1}...x_{r}e_{1}...e_{k}[a_{1},b_{1}]...[a_{h},b_{h}]$ or $x_{1}%
...x_{r}e_{1}...e_{k}d_{1}^{2}...d_{h}^{2}$, according as the orbit space is
orientable or not. Here $[a_{i},b_{i}]=a_{i}b_{i}a_{i}^{-1}b_{i}^{-1}$.

\section{Preliminary results}

Let $S$ be a (closed and orientable) surface of genus $g$, let $H_{1}%
(S,\mathbb{Z})$ be the first homology group of $S$ and $(.,.)$ be the bilinear
intersection form in $H_{1}(S,\mathbb{Z})$. Let $\sigma:S\rightarrow S$ be an
orientation reversing involution and assume that $Fix(\sigma)$ consists in $k$
ovals, $k\geq0$. We denote by $\sigma_{H}:H_{1}(S,\mathbb{Z})\rightarrow
H_{1}(S,\mathbb{Z})$ the homomorphism induced by $\sigma$ in homology.

\begin{lemma}
If $Fix(\sigma)$ is \textit{non-separating, i. e. }$S/\left\langle
\sigma\right\rangle $ is nonorientable, then there is a base $\left\{
X_{i},Y_{i}\right\}  _{i=1,...,g}$ in $H_{1}(S,\mathbb{Z})$ such that
$(X_{i},X_{j})=(Y_{i},Y_{j})=0$, $(X_{i},Y_{j})=\delta_{ij}$, the cycles
$X_{i}$, $i=1,...,k$ are represented by the ovals of $\sigma$, $\sigma
_{H}(X_{i})=X_{i}$, $i=1,...,g,$ and

$\sigma_{H}(Y_{i})=\left\{
\begin{array}
[c]{cc}%
-Y_{i}-\overset{g}{\underset{j=1}{\sum}}X_{j} & i=1,...,k\\
-Y_{i}-X_{i}-\overset{g}{\underset{j=1}{\sum}}X_{j} & i=k+1,...,g
\end{array}
\right.  $.
\end{lemma}

\textbf{Proof}. The orbifold $S/\left\langle \sigma\right\rangle $ can be
uniformized by a spherical, euclidean or hyperbolic crystallographic group
$\Delta$ of signature

\begin{center}
$(g-k+1,-,[-];\{(-),\overset{k}{...},(-)\})$,
\end{center}

\hspace{-1cm}i. e. $S/\left\langle \sigma\right\rangle =\mathbb{U}^{2}/\Delta
$. Let

\begin{center}
$\left\langle e_{1},...,e_{k},c_{1},...,c_{k},d_{k+1},...,d_{g+1}%
:e_{1}...e_{k}d_{k+1}^{2}...d_{g+1}^{2}=1,\right.  $

$\left.  c_{i}^{2}=1,e_{i}c_{i}e_{i}^{-1}=c_{i},i=1,...,k\right\rangle $
\end{center}

\hspace{-0.5cm}be a canonical presentation of $\Delta$.

The index two subgroup of $\Delta^{+}$ containing the orientation
transformations of $\Delta$, uniformizes $S$, i. e. $S=\mathbb{U}^{2}%
/\Delta^{+}$. A set of generators of $\Delta^{+}$ in terms of the above
presentation of $\Delta$ is:

\begin{center}
$\{e_{1},d_{g+1}c_{1},...,e_{k},d_{g+1}c_{k},d_{k+1}^{2},d_{g+1}%
d_{k+1},...,d_{g}^{2},d_{g+1}d_{g}\} $
\end{center}

The involution $\sigma$ induces on $\Delta^{+}$ an automorphism given by
conjugation by $d_{g+1}$. On the set of generators of $\Delta^{+}$ produces:

\begin{center}%
\begin{tabular}
[c]{ccccc}%
$e_{i}$ & $\rightarrow$ & $d_{g+1}e_{i}d_{g+1}^{-1}$ & $=$ & $(d_{g+1}%
c_{i})e_{i}(d_{g+1}c_{i})^{-1}$\\
$d_{i}^{2}$ & $\rightarrow$ & $d_{g+1}d_{i}^{2}d_{g+1}^{-1}$ & $=$ &
$(d_{g+1}d_{i})d_{i}^{2}(d_{g+1}d_{i})^{-1}$\\
$d_{g+1}c_{i}$ & $\rightarrow$ & $d_{g+1}^{2}c_{i}d_{g+1}^{-1}$ & $=$ &
$d_{g+1}^{2}(d_{g+1}c_{i})^{-1}$\\
$d_{g+1}d_{i}$ & $\rightarrow$ & $d_{g+1}^{2}d_{i}d_{g+1}^{-1}$ & $=$ &
$d_{g+1}^{2}d_{i}^{2}(d_{g+1}d_{i})^{-1}$%
\end{tabular}
\end{center}

Let $\theta:\Delta^{+}\rightarrow H_{1}(S,\mathbb{Z})$ be the abelianization
epimorphism and denote:

\begin{center}
$X_{i}=\theta(e_{i})$, $i=1,...,k$

$X_{i+k}=\theta(d_{i+k}^{2})$, $i=1,...,g-k$

$Y_{i}=\theta(d_{g+1}c_{i})$, $i=1,...,k$

$Y_{i+k}=\theta(d_{g+1}d_{i+k})$, $i=1,...,g-k$.
\end{center}

Considering the elements $e_{i}$ and $d_{i}$, $c_{i}$ in $\pi_{1}%
O(\mathbb{U}^{2}/\Delta)$, the above formula tell us how to repesent the
homology cycles $X_{i}$ and $Y_{i}$ as lifting of paths in $S/\left\langle
\sigma\right\rangle $ and then to compute the intersection form. We have:

\begin{center}
$(X_{i},X_{j})=(Y_{i},Y_{j})=0$, and $(X_{i},Y_{j})=\delta_{ij}$, for
$i,j\in\{1,...,g\}$.
\end{center}

The effect of $\sigma_{H}$ on $X_{i}$ is:

$\sigma_{H}(X_{i})=X_{i}$ and since $\theta(d_{g+1}^{2})=-\overset
{g}{\underset{i=1}{\sum}}X_{i}$ then:

\begin{center}
$\sigma_{H}(Y_{i})=\left\{
\begin{array}
[c]{cc}%
-Y_{i}-\overset{g}{\underset{j=1}{\sum}}X_{j} & i=1,...,k\\
-Y_{i}-X_{i}-\overset{g}{\underset{j=1}{\sum}}X_{j} & i=k+1,...,g
\end{array}
\right.  .$ $\ \ \ \square$
\end{center}

\begin{corollary}
Let $Z\in H_{1}(S,\mathbb{Z})$ such that $\sigma_{H}(Z)=Z$, then $Z=\sum
\alpha_{i}X_{i}$.
\end{corollary}

\textbf{Proof}. Let $Z=\sum\alpha_{i}X_{i}+\sum\beta_{i}Y_{i}$. If $\sigma
_{H}(Z)=Z$ then $-\beta_{i}=(Z,X_{i})=(\sigma_{H}(Z),X_{i})=\beta_{i}$. Hence
$\beta_{i}=0$. $\ \ \ \ \square$

We denote $V_{0}=\left\langle X_{1},...,X_{g}\right\rangle =Fix(\sigma_{H})$ and

$Aut^{\sigma_{H}}(H_{1}(S,\mathbb{Z}))=\{h\in Aut(H_{1}(S,\mathbb{Z})):h$
preserves the symplectic form $(.,.)$ and $h\circ\sigma_{H}=\sigma_{H}\circ
h\}$.

\begin{corollary}
If $g\in Aut^{\sigma_{H}}(H_{1}(S,\mathbb{Z}))$ then $g(V_{0})=V_{0}$ and
there is a homomorphism $\varphi:Aut^{\sigma_{H}}(H_{1}(S,\mathbb{Z}%
))\rightarrow Aut(V_{0})$ defined by $\varphi(a)=\left.  a\right|  _{V_{0}}$.
\end{corollary}

\textbf{Proof}. It is a consequence of the Corollary 2, since $V_{0}%
=Fix(\sigma_{H}).$ $\ \ \ \ \square$

\begin{corollary}
$\varphi$ is a monomorphism.
\end{corollary}

\textbf{Proof}. Let $a\in\ker\varphi$. Hence $a\in Aut^{\sigma_{H}}%
(H_{1}(S,\mathbb{Z}))$ and $\left.  a\right|  _{V_{0}}=id$.

We know that $a(X_{i})=X_{i}$, $i=1,...,g$ and now we want to show that
$a(Y_{i})=Y_{i}$.

Since $a$ preserves the symplectic form and $(X_{i},Y_{j})=\delta_{ij}$ then
$\delta_{ij}=(X_{i},Y_{j})=(a(X_{i}),a(Y_{j}))=(X_{i},a(Y_{j}))$, thus
$a(Y_{i})=Y_{i}+\sum\alpha_{ij}X_{j}$.

Now we use that $a\circ\sigma_{H}=\sigma_{H}\circ a$:

\begin{center}
$\sigma_{H}\circ a(Y_{i})=\sigma_{H}(Y_{i}+\sum\alpha_{ij}X_{j})=-Y_{i}%
-\varepsilon_{i}X_{i}-\overset{g}{\underset{i=1}{\sum}}X_{i}+\sum\alpha
_{ij}X_{j}$,
\end{center}

\hspace{-1cm}where $\varepsilon_{i}$ is 1 or 0,

\begin{center}
$a\circ\sigma_{H}(Y_{i})=a\circ(-Y_{i}-\varepsilon_{i}X_{i}-\overset
{g}{\underset{i=1}{\sum}}X_{i})=-Y_{i}-\sum\alpha_{ij}X_{j}-\varepsilon
_{i}X_{i}-\overset{g}{\underset{i=1}{\sum}}X_{i}.$
\end{center}

Thus $\alpha_{ij}=0$ and $a(Y_{i})=Y_{i}$, so $a=id_{H_{1}(S,\mathbb{Z})}$.
\ \ \ $\square$

\begin{lemma}
If $Fix(\sigma)$ is \textit{separating, i. e. }$S/\left\langle \sigma
\right\rangle $ is orientable (must be $k>0$) then there is a base $\left\{
X_{i},Y_{i}\right\}  _{i=1,...,g}$ in $H_{1}(S,\mathbb{Z})$ such that
$(X_{i},X_{j})=(Y_{i},Y_{j})=0$, $(X_{i},Y_{j})=\delta_{ij}$, the cycles
$X_{i}$, $i=1,...,k-1$ are represented by $k-1$ ovals of $\sigma$ and:

$\sigma_{H}(X_{i})=\left\{
\begin{array}
[c]{cc}%
X_{i} & i=1,...,k-1\\
X_{i+\frac{g-k+1}{2}} & i=k,...,\frac{g+k-1}{2}\\
X_{i-\frac{g-k+1}{2}} & i=\frac{g+k+1}{2},...,g
\end{array}
\right.  $

$\sigma_{H}(Y_{i})=\left\{
\begin{array}
[c]{cc}%
-Y_{i} & i=1,...,k-1\\
-Y_{i+\frac{g-k+1}{2}} & i=k,...,\frac{g-k-1}{2}\\
-Y_{i-\frac{g-k+1}{2}} & i=\frac{g-k+1}{2},...,g
\end{array}
\right.  $.
\end{lemma}

\textbf{Proof}. The proof is similar to the proof of Lemma 1 but now using the
uniformization of $S/\left\langle \sigma\right\rangle $ by a crystallographic
group $\Delta$ with signature $(\frac{g-k+1}{2},+,[-],\{(-),\overset{k}%
{...},(-)\})$.

Let

\begin{center}
$<e_{1}...e_{k},,c_{1},...,c_{k},a_{1},b_{1},...,a_{\frac{g-k+1}{2}}%
,b_{\frac{g-k+1}{2}}:$

$e_{1}...e_{k}[a_{1},b_{1}]...[a_{\frac{g-k+1}{2}},b_{\frac{g-k+1}{2}}]=1,$
$c_{i}^{2}=1,e_{i}c_{i}e_{i}^{-1}=c_{i},i=1,...,k>$
\end{center}

\hspace{-0.5cm}be a canonical presentation of $\Delta$.

The index two subgroup of $\Delta$, $\Delta^{+}$, consisting of the
orientation preserving transformations of $\Delta$ uniformizes $S$. A set of
generators of $\Delta^{+}$ in terms of the above presentation of $\Delta$ is:

\begin{center}
$\{e_{1},c_{k}c_{1},...,e_{k-1},c_{k}c_{k-1},a_{1},b_{1},...,a_{\frac
{g-k+1}{2}},b_{\frac{g-k+1}{2}},$

$c_{k}a_{1}c_{k},c_{k}b_{1}^{-1}c_{k},...,c_{k}a_{\frac{g-k+1}{2}}c_{k}%
,c_{k}b_{\frac{g-k+1}{2}}^{-1}c_{k}\}$
\end{center}

The involution $\sigma$ induces on $\Delta^{+}$ an automorphism given by
conjugation by $c_{k}$. On the set of generators of $\Delta^{+}$ produces:

\begin{center}%
\begin{tabular}
[c]{ccccc}%
$e_{i}$ & $\rightarrow$ & $c_{k}e_{i}c_{k}$ & $=$ & $(c_{k}c_{i})e_{i}%
(c_{k}c_{i})^{-1}$\\
$c_{k}c_{i}$ & $\rightarrow$ & $c_{i}c_{k}$ & $=$ & $(c_{k}c_{i})^{-1}$\\
$a_{i},b_{i}$ & $\rightarrow$ & $c_{k}a_{i}c_{k}$ & $,$ & $c_{k}b_{i}c_{k}$\\
$c_{k}a_{i}c_{k},c_{k}b_{i}^{-1}c_{k}$ & $\rightarrow$ & $a_{i}$ & $,$ &
$b_{i}^{-1}$%
\end{tabular}
\end{center}

Let $\theta:\Delta^{+}\rightarrow H_{1}(S,\mathbb{Z})$ be the abelianization
epimorphism and denote:

\begin{center}
$X_{i}=\theta(e_{i})$, $i=1,...,k-1,$

$Y_{i}=\theta(c_{k}c_{i})$, $i=1,...,k-1,$

$X_{i}=\theta(a_{i})$, $i=k,...,\frac{g+k-1}{2}$, $Y_{i}=\theta(b_{i})$,
$i=k,...,\frac{g+k-1}{2}$,

$X_{i}=\theta(c_{k}a_{i}c_{k})$, $i=\frac{g+k+1}{2},...,g$, $Y_{i}%
=\theta(c_{k}b_{i}^{-1}c_{k})$, $i=\frac{g+k+1}{2},...,g$.
\end{center}

Hence $(X_{i},X_{j})=(Y_{i},Y_{j})=0$ and $(X_{i},Y_{j})=\delta_{ij}$. The
effect of $\sigma_{H}$ on $X_{i}$ is as given in the statement of the Lemma.

$\hspace{12cm}\square$

Let $\Delta$ be a crystallographic\ group with signature

\begin{center}
$(\frac{g-k+1}{\eta},\pm,[-],\{(-),\overset{k}{...},(-)\})$,
\end{center}

\hspace{-0.5cm}where $\eta$ is 2 or 1 depending if there is a $+$ or a $-$ in
the signature, uniformizing $S/\left\langle \sigma\right\rangle $. Let
$\Delta^{+}$ be the index two subgroup of $\Delta$ containing the orientation
preserving transformations, then $\sigma:S=\mathbb{U}^{2}/\Delta
^{+}\rightarrow\mathbb{U}^{2}/\Delta^{+}=S$ and there is a $2$-fold orbifold
covering $S=\mathbb{U}^{2}/\Delta^{+}\rightarrow\mathbb{U}^{2}/\Delta
=S/\left\langle \sigma\right\rangle $.

\begin{lemma}
There is an isomorphism $\psi:Aut(\Delta)\rightarrow Mod_{g}^{\sigma}.$
\end{lemma}

\textbf{Proof}. Let $\alpha:\Delta\rightarrow\Delta$ be an automorphism, there
is a homeomorphism $h:\mathbb{U}^{2}/\Delta\rightarrow\mathbb{U}^{2}/\Delta$
inducing $\alpha$ (see \cite{MS} Corollary 7.4 for hyperbolic crystallographic
groups, for spherical and euclidean groups is a classical fact).

Let $h^{+}:\mathbb{U}^{2}/\Delta^{+}\rightarrow\mathbb{U}^{2}/\Delta^{+}$ be
the orientation preserving lifting of $h$ to $\mathbb{U}^{2}/\Delta^{+}$. Thus
$h^{+}\circ\alpha=\alpha\circ h^{+}$ and we define $\psi(\alpha)=h^{+}$.

Conversely, if $h^{+}:S=\mathbb{U}^{2}/\Delta^{+}\rightarrow\mathbb{U}%
^{2}/\Delta^{+}$ and $h^{+}\circ\alpha=\alpha\circ h^{+}$ then there is
$h:S/\left\langle \sigma\right\rangle =\mathbb{U}^{2}/\Delta\rightarrow
\mathbb{U}^{2}/\Delta$. The homomorphism $h$ is an orbifold automorphism that
produces an automorphism $\alpha$ of the orbifold fundamental group $\pi
_{1}O(S/\left\langle \sigma\right\rangle )=\Delta$.

$\hspace{12cm}\square$

\section{Poincar\'{e}'s theorem for real algebraic curves with empty real part.}

Let $S$ be a (closed and orientable) surface of genus $g$. Let $\sigma
:S\rightarrow S$ be an orientation reversing involution and assume that
$Fix(\sigma)=\emptyset$. We denote by $\sigma_{H}:H_{1}(S,\mathbb{Z}%
)\rightarrow H_{1}(S,\mathbb{Z})$ the homomorphism induced by $\sigma$ in homology.

By Lemma 1 there is a base $\left\{  X_{i},Y_{i}\right\}  _{i=1,...,g}$ in
$H_{1}(S,\mathbb{Z})$ such that $(X_{i},X_{j})=(Y_{i},Y_{j})=0$, $(X_{i}%
,Y_{j})=\delta_{ij}$ and

$\sigma_{H}(X_{i})=X_{i}$, $\sigma_{H}(Y_{i})=%
\begin{array}
[c]{cc}%
-Y_{i}-X_{i}-\overset{g}{\underset{j=1}{\sum}}X_{j}, & i=1,...,g
\end{array}
$.

By Lemma 2 if $Z\in H_{1}(S,\mathbb{Z})$ such that $\sigma_{H}(Z)=Z$, then
$Z=\sum\alpha_{i}X_{i}$. Hence $V_{0}=\left\langle X_{1},...,X_{g}%
\right\rangle =Fix(\sigma_{H})$ and we denote

$Aut^{\sigma_{H},0}(H_{1}(S,\mathbb{Z}))=\{h\in Aut(H_{1}(S,\mathbb{Z})):h$
preserves the symplectic form $(.,.)$ and $h\circ\sigma_{H}=\sigma_{H}\circ
h\}$.

By Corollaries 3 and 4 there is a monomorphism $\varphi:Aut^{\sigma_{H}%
,0}(H_{1}(S,\mathbb{Z}))\rightarrow Aut(V_{0})$ defined by $\varphi(a)=\left.
a\right|  _{V_{0}}$.

Let $H:Mod_{g}^{\sigma}\rightarrow Aut^{\sigma_{H},0}(H_{1}(S,\mathbb{Z}))$ be
the homomorphism given by $H(h)=h_{H}$, where $h_{H}$ is the homomorphism
induced in the homology by $h$, then our main result in this Section is:

\begin{theorem}
If $\sigma:S\rightarrow S$ is an orientation reversing involution such that
$Fix(\sigma)=\emptyset$ then $H(Mod_{g}^{\sigma})=Aut^{\sigma_{H},0}%
(H_{1}(S,\mathbb{Z}))$ is an epimorphism.
\end{theorem}

\textbf{Proof}. Let $\Delta$ be a crystallographic\ group with signature
$(g+1,-,[-],\{-\})$ uniformizing $S/\left\langle \sigma\right\rangle $. We
have $Aut(\Delta)\overset{\psi}{\rightarrow}Mod_{g}^{\sigma}\overset
{H}{\rightarrow}Aut^{\sigma_{H},0}(H_{1}(S,\mathbb{Z}))\overset{\varphi
}{\rightarrow}Aut(V_{0})$. By Corollary 4, $\varphi$ is a monomorphism then we
need only to prove that $\varphi\circ H\circ\psi$ is an epimorphism. If
$g\leq1$ the result is obvious since $V_{0}$ in such cases has dimension
$\leq1$. Let now assume that $g\geq2$.

To prove the Theorem we need to show that $\varphi\circ H\circ\psi$ is an epimorphism.

Let $\left\langle d_{1},...,d_{g+1}:d_{1}^{2}...d_{g+1}^{2}=1\right\rangle $
be a canonical presentation of $\Delta$. We define the automorphisms
$\alpha_{i}$, $i=1,...,g$, $\beta$ and $\gamma$ by:

\begin{center}
$\alpha_{i}:\left\{
\begin{tabular}
[c]{cccc}%
$d_{i}$ & $\rightarrow$ & $d_{i}^{2}d_{i+1}d_{i}^{-2}$ & \\
$d_{i+1}$ & $\rightarrow$ & $d_{i}$ & \\
$d_{j}$ & $\rightarrow$ & $d_{j}$ & $j\neq i,i+1$%
\end{tabular}
\right.  $

$\beta:\left\{
\begin{tabular}
[c]{cccc}%
$d_{1}$ & $\rightarrow$ & $d_{g+1}^{-1}$ & \\
$d_{2}$ & $\rightarrow$ & $d_{g}^{-1}$ & \\
$d_{j}$ & $\rightarrow$ & $d_{g-j+2}^{-1}$ & $j\geq3$%
\end{tabular}
\right.  $

$\gamma:\left\{
\begin{tabular}
[c]{cccc}%
$d_{1}$ & $\rightarrow$ & $d_{1}d_{2}^{-1}d_{1}^{-1}$ & \\
$d_{2}$ & $\rightarrow$ & $d_{1}d_{2}^{2}$ & \\
$d_{j}$ & $\rightarrow$ & $d_{j}$ & $j\geq3$%
\end{tabular}
\right.  $
\end{center}

Now we shall compute $\varphi\circ H\circ\psi$ of the above automorphisms. Let
$\left\{  X_{i}\right\}  _{i=1,...,g}$ be the basis of $V_{0}$ obtained from
the above canonical presentation of $\Delta$ following the proof of Lemma 1$.$
If $A_{i}=\varphi\circ H\circ\psi(\alpha_{i})$ and $X_{i}=\theta(d_{i}^{2})$
then $A_{i}(X_{i})=X_{i+1}$, $A_{i}(X_{i+1})=X_{i}$, $A_{i}(X_{j})=X_{j}$,
$j\neq i,i+1$.

The matrix of $A_{i}$ in the basis $\left\{  X_{i}\right\}  _{i=1,...,g}$ is:

\begin{center}%
\begin{tabular}
[c]{c}%
\\
\\
$i$\\
$i+1$%
\end{tabular}
$\left[
\begin{tabular}
[c]{cccccccccc}%
$1$ & $0$ & $...$ & $0$ &  &  &  &  &  & \\
$0$ & $1$ &  &  &  &  &  &  &  & \\
& $...$ & $...$ &  &  &  &  &  &  & \\
$0$ & $...$ & $0$ & $1$ &  &  &  &  &  & \\
&  &  &  & $0$ & $1$ &  &  &  & \\
&  &  &  & $1$ & $0$ &  &  &  & \\
&  &  &  &  &  & $1$ & $0$ & $...$ & $0$\\
&  &  &  &  &  & $0$ & $1$ &  & \\
&  &  &  &  &  &  & $...$ & $...$ & \\
&  &  &  &  &  & $0$ & $...$ & $0$ & $1$%
\end{tabular}
\right]  $,
\end{center}

If $B=\varphi\circ H\circ\psi(\beta)$, then $B(X_{1})=\overset{g}%
{\underset{i=1}{\sum}}X_{i}$ and $B(X_{j})=-X_{g-j+2}$. The matrix of $B$ in
the basis $\left\{  X_{i}\right\}  _{i=1,...,g}$ is:

\begin{center}
$\left[
\begin{tabular}
[c]{cccccc}%
$1$ & $0$ &  & $...$ & $0$ & $0$\\
$1$ & $0$ &  &  & $0$ & $-1$\\
$1$ & $0$ &  &  & $-1$ & $0$\\
&  &  & $...$ &  & \\
$1$ & $0$ & $-1$ &  & $0$ & $0$\\
$1$ & $-1$ & $0$ & $...$ & $0$ & $0$%
\end{tabular}
\right]  $.
\end{center}

And $C=\varphi\circ H\circ\psi(\gamma)$, $C(X_{1})=-X_{2}$, $C(X_{2}%
)=2X_{2}+X_{1}$, $C(X_{j})=X_{j}$, $j\neq1,2,$ with matrix in the basis
$\left\{  X_{i}\right\}  _{i=1,...,g}$:

\begin{center}
$\left[
\begin{tabular}
[c]{cccccc}%
$0$ & $1$ & $0$ & $...$ &  & $0$\\
$-1$ & $2$ & $0$ &  &  & $0$\\
$0$ & $0$ & $1$ & $0$ & $...$ & $0$\\
$0$ & $0$ & $0$ & $1$ &  & $0$\\
&  &  &  & $...$ & \\
$0$ & $0$ & $0$ & $...$ & $0$ & $1$%
\end{tabular}
\right]  $.
\end{center}

Note that multiplying by the left or the right by $A_{i}$ we can interchange
the rows and columns of a given matrix. The following matrix obtained from $B$
permuting rows corresponds to an element in the image of $\varphi\circ
H\circ\psi$:

\begin{center}
$B^{\prime}=\left[
\begin{tabular}
[c]{cccccc}%
$1$ & $0$ & $0$ & $...$ &  & $0$\\
$1$ & $-1$ & $0$ &  &  & $0$\\
$1$ & $0$ & $-1$ & $0$ & $...$ & $0$\\
$1$ & $0$ & $0$ & $-1$ &  & $0$\\
&  &  &  & $...$ & \\
$1$ & $0$ & $0$ & $...$ & $0$ & $-1$%
\end{tabular}
\right]  $
\end{center}

Now the following matrices also corresponds to elements in the image of
$\varphi\circ H\circ\psi$:

\begin{center}
$G_{1}=CA_{1}B^{\prime}=\left[
\begin{tabular}
[c]{cccccc}%
$1$ & $0$ & $0$ & $...$ &  & $0$\\
$1$ & $1$ & $0$ &  &  & $0$\\
$0$ & $0$ & $-1$ & $0$ & $...$ & $0$\\
$0$ & $0$ & $0$ & $-1$ &  & $0$\\
&  &  &  & $...$ & \\
$0$ & $0$ & $0$ & $...$ & $0$ & $-1$%
\end{tabular}
\right]  $
\end{center}

\hspace{-0.5cm}and:

\begin{center}
$G_{2}=\left[
\begin{tabular}
[c]{cccccc}%
$1$ & $1$ & $0$ & $...$ &  & $0$\\
$0$ & $1$ & $0$ &  &  & $0$\\
$0$ & $0$ & $-1$ & $0$ & $...$ & $0$\\
$0$ & $0$ & $0$ & $-1$ &  & $0$\\
&  &  &  & $...$ & \\
$0$ & $0$ & $0$ & $...$ & $0$ & $-1$%
\end{tabular}
\right]  .$
\end{center}

We have also:

$G_{3}=B^{\prime}CA_{1}B^{\prime}$ with matrix:

\begin{center}
$\left[
\begin{tabular}
[c]{cccccc}%
$1$ & $0$ & $0$ & $...$ &  & $0$\\
$0$ & $-1$ & $0$ &  &  & $0$\\
$0$ & $0$ & $1$ & $0$ & $...$ & $0$\\
$0$ & $0$ & $0$ & $1$ &  & $0$\\
&  &  &  & $...$ & \\
$0$ & $0$ & $0$ & $...$ & $0$ & $1$%
\end{tabular}
\right]  $.
\end{center}

Using $A_{i}$ and $G_{3}$ we can obtain from $G_{2}$ and $G_{3}$ the matrices:

\begin{center}
$G_{1}^{\prime}=\left[
\begin{tabular}
[c]{cccccc}%
$1$ & $0$ & $0$ & $...$ &  & $0$\\
$1$ & $1$ & $0$ &  &  & $0$\\
$0$ & $0$ & $1$ & $0$ & $...$ & $0$\\
$0$ & $0$ & $0$ & $1$ &  & $0$\\
&  &  &  & $...$ & \\
$0$ & $0$ & $0$ & $...$ & $0$ & $1$%
\end{tabular}
\right]  $ and $G_{2}^{\prime}=\left[
\begin{tabular}
[c]{cccccc}%
$1$ & $1$ & $0$ & $...$ &  & $0$\\
$0$ & $1$ & $0$ &  &  & $0$\\
$0$ & $0$ & $1$ & $0$ & $...$ & $0$\\
$0$ & $0$ & $0$ & $1$ &  & $0$\\
&  &  &  & $...$ & \\
$0$ & $0$ & $0$ & $...$ & $0$ & $1$%
\end{tabular}
\right]  $
\end{center}

To finish the proof it is enough to show that given an integer $g\times g$
matrix $M$ with $\det M=\pm1$ then $M$ is a product of the matrices $A_{i}$,
$G_{1}^{\prime}$ and $G_{2}^{\prime}$. Equivalently:

\begin{center}
$w_{1}(A_{i},G_{1}^{\prime},G_{2}^{\prime})Mw_{2}(A_{i},G_{1}^{\prime}%
,G_{2}^{\prime})=I$
\end{center}

\hspace{-0.5cm}where $w_{1}(A_{i},G_{1}^{\prime},G_{2}^{\prime})$ and
$w_{2}(A_{i},G_{1}^{\prime},G_{2}^{\prime})$ are words in $A_{i},B,C$ and this
fact follows from \cite{V}, vol. 2 pg. 107.

\section{Poincar\'{e}'s theorem for real Riemann surfaces with non-separating
fixed point set.}

Assume that $\sigma$ is an orientation reversing autohomeomorphism that fixes
$k$ closed curves and such that $S/\left\langle \sigma\right\rangle $ is
non-orientable. If $h$ is an autohomeomorphism of $S$ such that $h\circ
\sigma=\sigma\circ h$, then induces a homeomorphism $h^{\prime}:S/\left\langle
\sigma\right\rangle \rightarrow S/\left\langle \sigma\right\rangle $. Hence
$h^{\prime}$ sends boundary components to boundary components and $h$ sends
fixed curves by $\sigma$ in fixed curves, may be changing the orientation. Let
$X_{1},...,X_{k}$ be the cycles in $H_{1}(S,\mathbb{Z})$ represented by the
ovals of $\sigma$. Thus $h_{H}$ must permute the elements in the set $\left\{
\pm X_{1},...,\pm X_{k}\right\}  $, i. e. $h(X_{i})=\pm X_{j}$, $i,j\in
\{1,...,k\}$. Let $\left\{  X_{i},Y_{i}\right\}  _{i=1,...,g}$ be a base of
$H_{1}(S,\mathbb{Z})$ given by the Lemma 1. We shall denote:

\begin{center}
$Aut^{\sigma_{H},-k}(H_{1}(S,\mathbb{Z}))=\{h\in Aut(H_{1}(S,\mathbb{Z})):h$
preserves the symplectic form $(.,.)$, $h\circ\sigma_{H}=\sigma_{H}\circ h$,
$h(X_{i})=\pm X_{j}$, $i,j\in\{1,...,k\}$, $h(X_{i})=\overset{k}%
{\underset{j=1}{\sum}}2\alpha_{ij}X_{j}+\overset{g}{\underset{l=k+1}{\sum}%
}\beta_{il}X_{l}$, $i\in\{k+1,...,g\},\alpha_{ij},\beta_{il}\in\mathbb{Z}\}$

and

$Aut^{-k}(V_{0})=\{h\in Aut(V_{0}):h(X_{i})=\pm X_{j}$, $i,j\in\{1,...,k\},$

$h(X_{i})=\sum2\alpha_{ij}X_{j}+\sum\beta_{il}X_{l}$, $i\in\{k+1,...,g\}$,
$\alpha_{ij},\beta_{il}\in\mathbb{Z\}}$.
\end{center}

\bigskip

If $H:Mod_{g}^{\sigma}\rightarrow Aut(H_{1}(S,Z))$ is the homomorphism given
by $H(h)=h_{H}$, where $h_{H}$ is the homomorphism induced in the homology by
$h$ then $H(Mod_{g}^{\sigma})\subset Aut^{\sigma_{H},-k}(H_{1}(S,\mathbb{Z}%
))$. Following the notation in Section 3, note that $X_{i}=\theta(d_{i}d_{i})$
with $i\in\{k+1,...,g\}$. If $H(h)=h_{H}$, then $h_{H}(X_{i})=\theta(h_{\ast
}(d_{i})h_{\ast}(d_{i}))$, where $h_{\ast}$ is the homomorphism induced by $h$
on $\Delta$. Now $h_{\ast}(d_{i})$ is a word in the $d_{i}$ and the $e_{i}$
and then $\theta(h_{\ast}(d_{i})h_{\ast}(d_{i}))$ is a linear combination of
$2X_{j}=\theta(e_{j}^{2})$, $j\in\{1,...,k\}$ and $X_{l}=\theta(d_{l}^{2})$,
$l\in\{k+1,...,g\}$. Hence $H(Mod_{g}^{\sigma})\subset Aut^{\sigma_{H}%
,-k}(H_{1}(S,\mathbb{Z}))$.

The main result of this Section is:

\begin{theorem}
If $\sigma$ is an orientation reversing autohomeomorphism of a surface $S$
such that $\sigma$ fixes $k$ closed curves and $S/\left\langle \sigma
\right\rangle $ is non-orientable then $H(Mod_{g}^{\sigma})=Aut^{\sigma
_{H},-k}(H_{1}(S,\mathbb{Z}))$.
\end{theorem}

\textbf{Proof}. Assume that $S/\left\langle \sigma\right\rangle $ is
non-orientable and has $k\neq0$ boundary components.

Let $\left\{  X_{1},Y_{1},...,X_{g},Y_{g}\right\}  $ be a basis of
$H_{1}(S,\mathbb{Z)}$ given by Lemma 1 and let $h_{H}\in Aut^{\sigma_{H}%
,-k}(H_{1}(S,\mathbb{Z}))$.

\bigskip

\begin{lemma}
There are $h_{H}^{\prime}\in Aut^{\sigma_{H},-k}(H_{1}(S,\mathbb{Z}))$ and
$f\in Mod_{g}^{\sigma}$, such that $H(f)\circ h_{H}=h_{H}^{\prime}$ and
$h_{H}^{\prime}(X_{i})=X_{i}$, $i=1,...,k$ and the subspace $\left\langle
X_{k+1},...,X_{g}\right\rangle $ is invariant by the automorphism
$h_{H}^{\prime}.$
\end{lemma}

\textbf{Proof of the Lemma}. We define the automorphisms $\delta$, $\rho_{i}$
and $\mu_{i}$ of $\Delta$ by:

\begin{center}
$\delta:\left\{
\begin{tabular}
[c]{cccc}%
$e_{k}$ & $\rightarrow$ & $e_{k}d_{k+1}e_{k}^{-1}d_{k+1}^{-1}e_{k}^{-1}$ & \\
$d_{k+1}$ & $\rightarrow$ & $e_{k}d_{k+1}$ & \\
$c_{k}$ & $\rightarrow$ & $e_{k}d_{k+1}c_{k}d_{k+1}^{-1}e_{k}^{-1}$ & \\
$y$ & $\rightarrow$ & $y$ & for every canonical generator\\
&  &  & different from $e_{k}$, $d_{k+1}$ and $c_{k}$%
\end{tabular}
\right.  $

$\rho_{i}:\left\{
\begin{tabular}
[c]{cccc}%
$e_{i}$ & $\rightarrow$ & $e_{i}e_{i+1}e_{i}$ & \\
$e_{i+1}$ & $\rightarrow$ & $e_{i}$ & \\
$c_{i}$ & $\rightarrow$ & $e_{i}c_{i+1}e_{i}^{-1}$ & \\
$y$ & $\rightarrow$ & $y$ & for every canonical generator\\
&  &  & different from $e_{i}$, $e_{i+1}$ and $c_{i}$%
\end{tabular}
\right.  $

$\mu_{i}:\left\{
\begin{tabular}
[c]{cccc}%
$d_{i}$ & $\rightarrow$ & $d_{i}^{2}d_{i+1}d_{i}^{-2}$ & \\
$d_{i+1}$ & $\rightarrow$ & $d_{i}$ & \\
$y$ & $\rightarrow$ & $y$ & for every canonical generator\\
&  &  & different from $d_{i}$ and $d_{i+1}$%
\end{tabular}
\right.  $
\end{center}

The above automorphisms produce $D=\varphi\circ H\circ\varphi(\delta)$,
$R_{i}=\varphi\circ H\circ\varphi(\rho_{i})$, $M_{i}=\varphi\circ
H\circ\varphi(\mu_{i})$. In the basis $\{X_{l}\}_{l=1,...,g}$ of $V_{0}$ the
automorphisms $D,$ $R_{i}$ and $M_{i}$ give:

$D(X_{k})=-X_{k}$

$D(X_{k+1})=X_{k+1}+2X_{k}$

$D(X_{l})=X_{l}$, $l\neq k,k+1$,

\bigskip

$R_{i}(X_{i})=X_{i+1}$

$R_{i}(X_{i+1})=X_{i}$

$R_{i}(X_{l})=X_{l}$, $l\neq i,i+1$, $i\in\{1,...,k-1\},$

$\bigskip$

$M_{i}(X_{i})=X_{i+1}$

$M_{i}(X_{i+1})=X_{i}$

$M_{i}(X_{l})=X_{l}$, $l\neq i,i+1$, $i\in\{k+1,...,g-1\}.$

\bigskip

We can construct a word $w(D,R_{i},M_{i})$ in $D,R_{i},M_{i}$ such that:

\begin{center}
$w(D,R_{i},M_{i})(h(X_{j}))=X_{j}$, $j=1,...,k$ and

$w(D,R_{i},M_{i})(h(X_{l}))\in\left\langle X_{k+1},...,X_{g}\right\rangle $,
$l\in\{k+1,...,g\}$.
\end{center}

Hence the element $f\in Mod_{g}^{\sigma}$ that we are looking for is
$\psi(w(\delta,\rho_{i},\mu_{i}))$.

$\hspace{10cm}\square$

\bigskip

We shall denote by $h_{H,1}$ the restriction of $h_{H}^{\prime}$ to
$\left\langle X_{k+1},...,X_{g}\right\rangle $.

Let $C$ be a closed curve in $S/\left\langle \sigma\right\rangle $ such that
$S/\left\langle \sigma\right\rangle -C$ has two connected components $F_{1}$
and $F_{2},$ such that $\overline{F_{1}}$ is homeomorphic to a non-orientable
surface of genus $g+1$ with connected boundary and $\overline{F_{2}}$ is
homeomorphic to a planar surface with $k+1$ boundary components. Note that
$C=\overline{F_{1}}\cap\overline{F_{2}}$ and that the preimage by
$\pi:S\rightarrow S/\left\langle \sigma\right\rangle $ of $C$ is a set of two
closed and disjoint curves $C_{1}$ and $C_{2}$.

Let $\widehat{F_{1}}$ be the surface obtained from $\overline{F_{1}}$ closing
the boundary component with a disc and $\widetilde{F_{1}}$ be the surface
obtained from $\pi^{-1}(F_{1})$ closing with two discs the boundaries $C_{1}$
and $C_{2}$. Let $\sigma_{1}$ be the orientation reversing automorphism of
$\widetilde{F_{1}}$ giving as quotient $\widehat{F_{1}}$. Hence we can
identify $Fix(\sigma_{1})\subset H_{1}(\widetilde{F_{1}},\mathbb{Z})$ with
$\left\langle X_{k+1},...,X_{g}\right\rangle $. By Theorem 7 we have a
homeomorphism $h_{1}:\widetilde{F_{1}}\rightarrow\widetilde{F_{1}}$ such that
$h_{1}\circ\sigma_{1}=\sigma_{1}\circ h_{1}$, $(h_{1})_{H}=h_{H,1}$ and we can
choose $h_{1}$ in such a way that $h_{1}$ fixes $C_{1}$ and $C_{2}$. The
extension of $\left.  h_{1}\right|  _{\pi^{-1}(F_{1})}$ to $S$ by the identity
give us an element $h^{\prime}$ of $Mod_{g}^{\sigma}$ such that $H(h^{\prime
})=h_{H}^{\prime}$.

$\hspace{10cm}\square$

\textbf{Example}. Let $S$ be a surface of genus $2$ and $\sigma$ be an
orientation reversing autohomeomorphism of $S$ with one (non-separating) oval.
Let $\{X_{1},Y_{1},X_{2},Y_{2}\}$ be a base given by the Lemma 1. Hence:
$\sigma_{H}(X_{i})=X_{i}$ and

$\sigma_{H}(Y_{1})=-Y_{1}-X_{1}-X_{2}$, $\sigma_{H}(Y_{2})=-Y_{2}-X_{1}%
-2X_{2}$.

Let us consider the automorphism of $H_{1}(S,\mathbb{Z})$:

$h(X_{1})=X_{1}$, $h(X_{2})=X_{1}+X_{2}$

$h(Y_{1})=-X_{1}-X_{2}+Y_{1}-Y_{2}$, $h(Y_{2})=Y_{2}-X_{1}$.

It is easy to verify that $h$ preserves the symplectic form $(.,.)$ and
$h\circ\sigma_{H}=\sigma_{H}\circ h$, but, by Theorem 8, $h$ is not induced by
an autohomeomorphism of $S$ that commutes with $\sigma$.

\section{Poincar\'{e}'s theorem for real Riemann surfaces with separating
fixed point set.}

Assume that $\sigma$ is an orientation reversing autohomeomorphism that fixes
$k$ closed curves and such that $S/\left\langle \sigma\right\rangle $ is
orientable, then $S-Fix(\sigma)$ has two connected components $S_{1}$ and
$S_{2}$. If $h$ is an orientation preserving autohomeomorphism of $S$ such
that $h\circ\sigma=\sigma\circ h$, $h$ induces an orientation preserving
homeomorphism $h^{\prime}:S_{i}\rightarrow S_{j}$ where $i=j$ or $i\neq j$
(depending on the nature of $h$). Hence $h^{\prime}$ sends boundary components
to boundary components and $h$ sends fixed curves by $\sigma$ in fixed curves.
Let $X_{1},...,X_{k-1},\overset{k-1}{\underset{i=1}{\sum}}X_{i}$ be the cycles
in $H_{1}(S,\mathbb{Z})$ represented by the ovals of $\sigma$. Thus $h_{H}$
satisfy $h(X_{i})=\varepsilon X_{j}$ or $\varepsilon\overset{k-1}%
{\underset{i=1}{\sum}}X_{i}$, $i,j\in\{1,...,k-1\}$, $\varepsilon=1$ if
$h^{\prime}:S_{i}\rightarrow S_{i}$ and $\varepsilon=-1$ if $h^{\prime}%
:S_{i}\rightarrow S_{j},$ $i\neq j$. We can identify

$H_{1}(S_{1},\mathbb{Z})=\left\langle X_{1},...,X_{k-1},X_{k},Y_{k}%
,...,X_{\frac{g-k-1}{2}},Y_{\frac{g-k-1}{2}}\right\rangle $ and

$H_{1}(S_{2},\mathbb{Z})=\left\langle X_{1},...,X_{k-1},X_{\frac{g-k+1}{2}%
},Y_{\frac{g-k+1}{2}},...,X_{g},Y_{g}\right\rangle $.

If we denote by

$H_{1}=\left\langle X_{k},Y_{k},...,X_{\frac{g-k-1}{2}},Y_{\frac{g-k-1}{2}%
}\right\rangle $ and $H_{2}=\left\langle X_{\frac{g-k+1}{2}},Y_{\frac
{g-k+1}{2}},...,X_{g},Y_{g},\right\rangle $,

\hspace{-0.5cm}then the homology classes in $H_{i}$ are represented by curves
contained in one of the connected components $S_{1}$ or $S_{2}$. Then
$h_{H}(H_{i})\subset H_{i}$ $i,j\in\{1,2\}$ if $\varepsilon=1$ and
$h_{H}(H_{i})\subset H_{j}$ $i\neq j$ if $\varepsilon=-1$. We shall denote:

\begin{center}
$Aut^{\sigma_{H},+k}(H_{1}(S,\mathbb{Z}))=\{h\in Aut(H_{1}(S,\mathbb{Z})):h$
preserves the symplectic form $(.,.)$, $h\circ\sigma_{H}=\sigma_{H}\circ h$,
$h(X_{i})=\varepsilon X_{j}$ or $\varepsilon\overset{k-1}{\underset{i=1}{\sum
}}X_{i},$ $i,j\in\{1,...,k-1\},$ $\varepsilon=\pm1$, $h(H_{i})\subset H_{i}$
if $\varepsilon=1$ and $h(H_{i})\subset H_{j},$ $i\neq j$, if $\varepsilon=-1$ $\},$
\end{center}

If $H:Mod_{g}^{\sigma}\rightarrow Aut(H_{1}(S,Z))$ is the homomorphism given
by $H(h)=h_{H}$, we have observed above that $H(Mod_{g}^{\sigma})\subset
Aut^{\sigma_{H},+k}(H_{1}(S,\mathbb{Z}))$.

\begin{theorem}
If $\sigma$ is an orientation reversing autohomeomorphism such that $\sigma$
fixes $k$ closed curves and such that $S/\left\langle \sigma\right\rangle $ is
orientable then $H(Mod_{g}^{\sigma})=Aut^{\sigma_{H},+k}(H_{1}(S,\mathbb{Z}))$.
\end{theorem}

\textbf{Proof}. Assume that $S/\left\langle \sigma\right\rangle $ is
orientable and has $k\neq0$ boundary components.

Let $\left\{  X_{1},Y_{1},...,X_{g},Y_{g}\right\}  $ be the basis of
$H_{1}(S,\mathbb{Z)}$ given in Lemma 4 and let $h_{H}\in Aut^{\sigma_{H}%
,+k}(H_{1}(S,\mathbb{Z}))$.

\bigskip

\begin{lemma}
There are $h_{H}^{\prime}\in Aut^{\sigma_{H},+k}(H_{1}(S,\mathbb{Z}))$ and
$f\in Mod_{g}^{\sigma}$, such that $H(f)\circ h_{H}=h_{H}^{\prime}$,
$h_{H}^{\prime}(X_{i})=X_{i}$, $i=1,...,k$ and $h_{H}^{\prime}(H_{i})=H_{j}$,
$i,j\in\{1,2\}$.
\end{lemma}

\textbf{Proof of the Lemma}. We define the automorphisms $\omega$, $\xi$,
$\nu_{i}$, $\mu$, $\rho_{i}$ of $\Delta$ by:

\begin{center}
$\omega:\left\{
\begin{tabular}
[c]{cccc}%
$a_{1}$ & $\rightarrow$ & $a_{1}b_{1}$ & \\
$y$ & $\rightarrow$ & $y$ & for every canonical generator\\
&  &  & different from $a_{1}$%
\end{tabular}
\right.  $

$\tau:\left\{
\begin{tabular}
[c]{cccc}%
$a_{1}$ & $\rightarrow$ & $a_{1}b_{1}$ & \\
$b_{1}$ & $\rightarrow$ & $a_{1}^{-1}$ & \\
$y$ & $\rightarrow$ & $y$ & for every canonical generator\\
&  &  & different from $a_{1}$, $b_{1}$%
\end{tabular}
\right.  $

$\nu_{i}:\left\{
\begin{tabular}
[c]{cccc}%
$a_{i}$ & $\rightarrow$ & $a_{i+1}$ & \\
$b_{i}$ & $\rightarrow$ & $b_{i+1}$ & \\
$a_{i+1}$ & $\rightarrow$ & $w_{i+1}^{-1}a_{i}w_{i+1}$ & where $w_{i+1}%
=[a_{i+1},b_{i+1}]$\\
$b_{i+1}$ & $\rightarrow$ & $w_{i+1}^{-1}b_{i}w_{i+1}$ & \\
$y$ & $\rightarrow$ & $y$ & for every canonical generator\\
&  &  & different from $a_{i}$, $a_{i+1}$, $b_{i}$, $b_{i+1}$%
\end{tabular}
\right.  $

$\mu:\left\{
\begin{tabular}
[c]{cccc}%
$e_{k}$ & $\rightarrow$ & $a_{1}^{-1}e_{k}a_{1}$ & \\
$a_{1}$ & $\rightarrow$ & $[a_{1}^{-1},e_{k}^{-1}]a_{1}$ & \\
$b_{1}$ & $\rightarrow$ & $b_{1}a_{1}^{-1}e_{k}a_{1}$ & \\
$y$ & $\rightarrow$ & $y$ & for every canonical generator\\
&  &  & different from $e_{k}$, $a_{1}$ and $b_{1}$%
\end{tabular}
\right.  $

$\rho_{i}:\left\{
\begin{tabular}
[c]{cccc}%
$e_{i}$ & $\rightarrow$ & $e_{i}e_{i+1}e_{i}$ & \\
$e_{i+1}$ & $\rightarrow$ & $e_{i}$ & \\
$c_{i}$ & $\rightarrow$ & $e_{i}c_{i+1}e_{i}^{-1}$ & \\
$c_{i+1}$ & $\rightarrow$ & $c_{i}$ & \\
$y$ & $\rightarrow$ & $y$ & for every canonical generator\\
&  &  & different from $e_{i}$, $e_{i+1}$ and $c_{i}$%
\end{tabular}
\right.  $
\end{center}

Having in consideration that $\theta(e_{k})=-\overset{k-1}{\underset{i=1}%
{\sum}}X_{i}$, the above automorphisms produce $Z=\varphi\circ H\circ
\varphi(\omega)$, $T=\varphi\circ H\circ\varphi(\tau)$, $N_{i}=\varphi\circ
H\circ\varphi(\nu_{i})$, $M=\varphi\circ H\circ\varphi(\mu)$, $R_{i}%
=\varphi\circ H\circ\varphi(\rho_{i})$ described by the formulae:

\bigskip

$Z(X_{k})=X_{k}+Y_{k},$ $Z(X_{k+\frac{g-k+1}{2}})=X_{k+\frac{g-k+1}{2}%
}-Y_{k+\frac{g-k+1}{2}}$

\bigskip

$T(X_{k})=X_{k}+Y_{k},$ $T(Y_{k})=-X_{k}$

$T(X_{k+\frac{g-k+1}{2}})=X_{k+\frac{g-k+1}{2}}-Y_{k+\frac{g-k+1}{2}},$
$T(Y_{k+\frac{g-k+1}{2}})=X_{k+\frac{g-k+1}{2}}$

\bigskip

$N_{i}(X_{i})=X_{i+1},$ $N_{i}(Y_{i})=Y_{i+1}$

$N_{i}(X_{i+1})=X_{i},$ $N_{i}(Y_{i+1})=Y_{i}$

$N_{i}(X_{i+\frac{g-k+1}{2}})=X_{i+1+\frac{g-k+1}{2}},$ $N_{i}(Y_{i+\frac
{g-k+1}{2}})=Y_{i+1+\frac{g-k+1}{2}}$

$N_{i}(X_{i+1+\frac{g-k+1}{2}})=X_{i+\frac{g-k+1}{2}},$ $N_{i}(Y_{i+1+\frac
{g-k+1}{2}})=Y_{i+\frac{g-k+1}{2}}$

\bigskip

$M(Y_{k})=Y_{k}-\overset{k-1}{\underset{i=1}{\sum}}X_{i},$ $M(Y_{k+\frac
{g-k+1}{2}})=Y_{k+\frac{g-k+1}{2}}-\overset{k-1}{\underset{i=1}{\sum}}X_{i}$

\bigskip

$R_{i}(X_{i})=X_{i+1},$ $R_{i}(Y_{i})=Y_{i+1}$

$R_{i}(X_{i+1})=X_{i},$ $R_{i}(Y_{i+1})=Y_{i}$

$i\in\{1,...,k-2\},$

\bigskip

$R_{k-1}(X_{k-1})=-\overset{k-1}{\underset{i=1}{\sum}}X_{i},$ $R_{k-1}%
(Y_{k-1})=-Y_{k-1}$

$R_{k-1}(Y_{i})=Y_{i}-Y_{k-1}$, $i=1,...,k-2$.

\bigskip

In the above formulae we have just written the image of the elements that
change by the automorphisms.

Finally we consider $U$ induced by an orientation preserving involution $u$ in
$S$ with $\sigma\circ u=u\circ\sigma$ such that $U(S_{1})=S_{2}$, $u$ has two
fixed points on each oval of $\sigma$ and

$U(X_{i})=-X_{i}$, $U(Y_{i})=-Y_{i}$, $i=1,...,k-1,$

$U(X_{i})=X_{i+\frac{g-k+1}{2}},$ $U(Y_{i})=Y_{i+\frac{g-k+1}{2}},$
$i=k,...,\frac{g+k-1}{2},$

$U(X_{i})=X_{i-\frac{g-k+1}{2}}$, $U(Y_{i})\rightarrow Y_{i-\frac{g-k+1}{2}}$,
$i=\frac{g+k+1}{2},...,g$.\bigskip

Now, as in the proof of Lemma 9, we can construct a word $w(Z,T,N_{i}%
,M,R_{i},U)$ in $Z,T,N_{i},M,R_{i},U$ such that:

\begin{center}
$w(Z,T,N_{i},M,R_{i},U)(h(X_{j}))=X_{j}$, $j=1,...,k-1$ and

$w(Z,T,N_{i},M,R_{i},U)(H_{i})\subset H_{i},$ $i,j\in\{1,2\}$.
\end{center}

Hence the element $f\in Mod_{g}^{\sigma}$ that we are looking for is
$\psi(w(\omega,\tau,\nu_{i},\mu,\rho_{i}))$.

$\hspace{10cm}\square$

\bigskip

We shall denote by $h_{H,1}$ the restriction of $h_{H}^{\prime}$ to $H_{1}$.
Let $C$ be a closed curve in $S_{1}$ such that $S_{1}-C$ has two connected
components $F_{1}$ and $F_{2},$ such that $\overline{F_{1}}$ is homeomorphic
to an orientable surface of genus $\frac{g-k+1}{2}$ with connected boundary
and $\overline{F_{2}}$ is homeomorphic to a planar surface with $k+1$ boundary components.

Let $\widehat{F_{1}}$ be the surface obtained from $\overline{F_{1}}$ closing
the boundary component with a disc. Hence we can identify $H_{1}%
(\widetilde{F_{1}},\mathbb{Z})$ with $H_{1}$. By Poincar\'{e} theorem there is
$h_{1}:\widetilde{F_{1}}\rightarrow\widetilde{F_{1}}$ such that $(h_{1}%
)_{H}=h_{H,1}$ and $h_{1}$ fixes $C_{1}$. We construct now

$h^{\prime}(z)=h_{1}(z)$ for $z\in\widetilde{F_{1}}\cap S_{1}$,

$h^{\prime}(z)=z$ for $z\in S_{1}-\widetilde{F_{1}},$

$h^{\prime}(z)=\sigma\circ h^{\prime}\circ\sigma$ $(z)$ for $z\in S_{2}$.

Thus we have an element $h^{\prime}$ of $Mod_{g}^{\sigma}$ such that
$H(h^{\prime})=h_{H}^{\prime}$.

$\hspace{10cm}\square$

{\small \bigskip}

{\small Antonio F. Costa}

{\small Departamento de Matem\'{a}ticas Fundamentales, Facultad de Ciencias,
UNED, C. Senda del rey, 9, 28040-Madrid, SPAIN}

{\small \bigskip}

{\small Sergey M. Natanzon}

{\small Mathematics Department, Independent University of Moscow, Belozersky
Institute, Moscow State University Moscow, RUSSIA}
\end{document}